\documentclass[12pt,twoside]{article}
\def\tilde{\widetilde}

\def\Tr{{\rm Tr}} 

\begin{document}
%\numberwithin{equation}{section}
%\def\ee{\end{equation}}
\title{Conformal Field Theory In Four And Six Dimensions}

\author{Edward Witten\thanks{Supported in part by NSF Grant
PHY-0070928. Lectures at the conference on  Topology,  Geometry, and Quantum
Field Theory, Oxford University (July, 2002).} \\
Institute for Advanced Study\\ Princeton, NJ 08540 USA}

\maketitle 
%%%%%%%%%%%%%%%%%%%%%%%%%%%%%%%%%%%%%%%%%%%%%%%%%%%%%%%%%%%%%%%%%%%%%

%\begin{abstract}
%No Abstract in this paper 
%\end{abstract}

\section{ Introduction}
%%%%%%%%%%%%%%%%%%%%%%%%%%%%%%%%%%%%%%%%%%%%%%%%

In these lectures, I will be considering  conformal field theory (CFT) mainly
in four and six dimensions, occasionally recalling facts about two
dimensions.  The notion of conformal field theory is familiar to
physicists.  From a mathematical point of view, we can keep in
mind Graeme Segal's definition \cite{gsegal} of conformal field theory.
Instead of just summarizing the definition here, I will review how physicists
actually study examples of quantum field theory, as this will make
clear the motivation for the definition.

When possible (and we will later consider examples in which this is
not possible), physicists make models of quantum field theory
using path integrals.  This means first of all
that, for any $n$-manifold $M_n$, we are given a space of fields  on
$M_n$; let us call the fields $\Phi$.  The fields might be, for example, real-valued functions, or
gauge fields (connections on a $G$-bundle over $M_n$ for some fixed
Lie group $G$), or $p$-forms on $M_n$ for some fixed $p$, or they
might be maps $\Phi:M_n\to W$ for some fixed manifold $W$. 
Then we are given a local  action functional $I(\Phi)$. ``Local'' means that the Euler-Lagrange
equations for a critical point of $I$ are partial differential equations. If we are
constructing a quantum field theory that is not required to be conformally
invariant, $I$ may be defined using a metric on $M_n$.  For conformal
field theory, $I$ should be defined using only a conformal structure.
For a closed $M_n$,
the  partition function $Z(M_n)$ is defined, formally, as the integral
over all $\Phi$ of $e^{-I(\Phi)}$:
\begin{equation}
Z(M_n)=\int D\Phi\,\exp(-I(\Phi)).
\label{phis}
\end{equation}
If $M_n$ has a boundary $M_{n-1}$, the integral depends on the boundary
conditions.  If we let
$\varphi$ denote the restriction of $\Phi$ to $M_{n-1}$, 
then it formally makes sense to consider
a  path integral on a manifold with boundary in which we integrate
over all $\Phi$ for some fixed $\varphi$. This defines
a function 
\begin{equation}
\Psi(\varphi)=\int_{\Phi|_{M_{n-1}}=\varphi} D\Phi\,\exp(-I(\Phi)).
\label{varphi}
\end{equation}
We interpret the function $\Psi(\varphi)$ as a vector in a Hilbert
space ${\cal H}(M_{n-1})$ of ${\bf L}^2$ functions of 
$\varphi$. 
{}From this starting point, one can motivate the sort of axioms for
quantum field theory that Segal considered.  I will not go into details,
as we will not need them in the present lectures.  In fact, to keep things
simple, we will mainly consider closed manifolds $M_n$ and the partition
function $Z(M_n)$.

Before getting to the specific examples that we will consider,
I will start with a general survey of conformal field theory in
various dimensions.
Two-dimensional conformal field theory plays an important role
in string theory and statistical mechanics and is also
relatively familiar mathematically.
\footnote{In counting dimensions, we include time, so a two-dimensional
theory, if formulated in Lorentz signature,
 is a theory in a world of one space and one time dimension.
In these lectures, we will mostly work with Euclidean signature.}
For example, rational conformal field theory is studied in detail
using complex geometry.  More general conformal field theories underlie,
for example, mirror symmetry.

Three and four-dimensional conformal field theory is also important for
physics.  Three-dimensional conformal field theory is used to describe
second order phase transitions in equilibrium statistical mechanics,
and a four-dimensional conformal field theory could conceivably play
a role in models of elementary particle physics.

Physicists used to think that four was the maximum dimension for
non-trivial (or non-Gaussian) unitary conformal field theory.  
Initially, therefore, little note was taken of a result by
Nahm \cite{wnahm} which implies that {\it six} 
is  the maximum possible dimension in the
supersymmetric case.  (A different result proved in the same paper
-- eleven is the maximal possible dimension for supergravity -- had a large
impact right away.)  Nahm's result follows from an algebraic argument
and I will explain what it says in section 3.  String theorists have
been quite surprised in the last few years to learn that the higher
dimensional superconformal field theories whose existence is suggested
by Nahm's theorem apparently do exist.  Explaining this, or at least giving
 a few hints, is the goal of these lectures.

One of the surprises is that the new theories suggested by Nahm's theorem
are theories for which there is apparently no Lagrangian -- at least none that
can be constructed using classical variables of any known sort.
Yet these new theories are intimately connected with fascinating mathematics
and physics of more conventional theories in four dimensions.

In section 2, we warm up with some conventional and less conventional
linear theories.  Starting with the example of abelian gauge theory in
four dimensions, I will describe some free or in a sense linear
conformal field theories
that can be constructed in arbitrary even dimensions.  The cases of
dimension $4k$ and $4k+2$ are rather different, as we will see.  The
most interesting linear theory in $4k+2$ dimensions is a self-dual theory
that does not have a Lagrangian, 
yet it exists quantum mechanically and its existence
is related to subtle modular behavior 
of the linear theories in $4k$ dimensions.

In section 3, I will focus on certain nonlinear examples in four and six
dimensions and the relations between them.  These examples will
be supersymmetric.  The importance for
us of supersymmetry is that it gives severe constraints that have 
made it possible to get
some insight about  highly nonlinear theories.   
After reviewing Nahm's theorem, I will
say a word or two about supersymmetric gauge theories in four dimensions
that are conformally invariant at the quantum level, and then about
how some of them are apparently related to nonlinear
superconformal field theories
in six dimensions.

\section{Gauge Theory And Its Higher Cousins}

First let us review abelian gauge theory, with gauge group $U(1)$.  
(For general references
on some of the following discussion of abelian gauge fields and
self-dual $p$-forms, see \cite{ewitten}.) The connection $A$ is locally
a one-form.  Under a gauge transformation, it transforms by
$A\to A+d\epsilon$, with $\epsilon$ a zero-form.  The curvature $F=dA$
is invariant.

  For the action, we take
\begin{equation}
I(A)={1\over 2e^2}\int_M F\wedge *F +{i\theta\over 2}\int_M
{F\over 2\pi}\wedge {F\over 2\pi}.
\label{wedge}
\end{equation}
Precisely in four dimensions, the Hodge $*$ operator on two-forms is
conformally invariant and so $I(A)$ is conformally invariant.  If $M$ is
closed, the  second
term in $I(A)$ is a topological invariant, $i(\theta/2)\int_M c_1({\cal L})^2.$
In general, $c_1({\cal L})^2$ is integral, and on a spin manifold
it is actually even.  So the integrand $\exp(-I(A))$ of the partition
function is always invariant to $\theta\to \theta+4\pi$, while
on a spin manifold it is invariant to $\theta\to\theta+2\pi$. 
In general, even when $M$ is not closed, this is a symmetry of the theory
(but in case $M$ has a boundary, the discussion becomes a little more
elaborate).

Now let us look at the partition function $Z(M)=\sum_{\cal L}\int DA\,\exp(-I(A))$, 
where we understand the sum over all possible connections $A$ as including
a sum over the line bundle ${\cal L}$ on which $A$ is a connection.  
We can describe the path integral rather explicitly, using the
decomposition $A=A'+A_h^{\cal L}$, where $A'$ is a connection on a trivial line
bundle ${\cal O}$, and $A_h^{\cal L}$ is (any) connection on ${\cal L}$
of harmonic curvature $F_h^{\cal L}$.  The action is
$I(A)=I(A')+I(A_h^{\cal L}),$
and the path integral is
\begin{equation}
\sum_{\cal L}\int DA\,\exp(-I(A))=\int DA'\,\exp(-I(A'))\sum_{\cal L}
\exp(-I(A_h^{\cal L})).
\label{facto}
\end{equation}
Here, note that $A_h^{\cal L}$ depends on ${\cal L}$, but $A'$ does not.

Let us look first at the second factor in eqn. \ref{facto}, the sum over ${\cal L}$.
On the lattice $H^2(M;{\bf Z})$, there is a natural, generally indefinite
quadratic form given, for $x$ an integral harmonic two-form, by
$(x,x)=\int_M x\wedge x$.  
There is also a positive-definite but
metric-dependent form $\langle x,x\rangle =\int_M x\wedge *x$,
with $*$ being the Hodge star operator.  The indefinite form $(x,x)$ has
signature $(b_{2,+},b_{2,-})$, where $b_{2,\pm}$ are the dimensions of the spaces
of self-dual and anti-self-dual harmonic two-forms.

Setting $x=F_h^{\cal L}/2\pi$, the sum over line bundles becomes
\begin{equation}
\sum_{x\in H^2(M;{\bf Z})} \exp\left(\begin{array}{c}-{4\pi^2\over e^2}
\langle x,x\rangle
+i{\theta\over 2}(x,x) \label{rang} \end{array} \right).
\end{equation}
If I set
$\tau ={\theta\over 2\pi}+{4\pi i\over e^2},$
then this function has modular properties with respect to $\tau$.
It is the non-holomorphic theta function of C. L. Siegel, which in the mid-1980's was
introduced in string theory by K. S. Narain to understand toroidal
compactification of the heterotic string.
The Siegel-Narain function has a simple transformation law under the
full modular group $SL(2,{\bf Z})$ if $M$ is spin, in which case
$(x,x)/2$ is integer-valued.  In general, it has modular properties for
a subgroup $\Gamma_0(2)$ of $SL(2,{\bf Z})$.  In any case, it transforms
as a modular function with holomorphic and anti-holomorphic weights
$(b_{2,+},b_{2,-})$.

The other factor in eqn. \ref{facto}, namely the integral over $A'$,
$\int DA' \,\exp(-I(A')),$
is essentially a Gaussian integral that can be defined by zeta functions.
Its dependence on the metric of $M$ is very complicated, but its dependence
on $\tau$ is very simple -- just a power of ${\rm Im}\,\tau$.  Including
this factor, the full path-integral transforms as a modular function
of weights $(1-b_1+b_{2,+}/2,1-b_1+b_{2,-}/2)
=((\chi+\sigma)/2,(\chi-\sigma)/2)$,
where $b_1$, $\chi$, and $\sigma$ are respectively the first Betti
number, the Euler characteristic, and the signature of $M$.

The
fact that the modular weights are linear combinations of $\chi$ and $\sigma$
has an important consequence, which I will not be able to explain fully
here.
Because $\chi$ and $\sigma$ can be written as integrals over $M$ of
quadratic polynomials in the Riemann curvature (using for example the
Gauss-Bonnet-Chern formula for $\chi$), it is possible to add
to the action $I$ a ``$c$-number'' term -- the integral of a local
expression that depends  on $\tau$
and on the metric of $M$ but not on the integration variable $A$ of the
path integral -- that cancels the modular weight and makes the partition
function completely invariant under $SL(2,{\bf Z})$ or $\Gamma_0(2)$.
The appropriate $c$-number terms arise naturally when, as we discuss
later, one derives the
four-dimensional abelian gauge theory from a six-dimensional self-dual
theory.

\bigskip\noindent{\it $p$-Form Analog}

Now I want to move on to the $p$-form analog,
for $p>2$.  For our purposes, we will be  informal in describing
$p$-form fields.
A ``$p$-form field'' $A_p$ is an object that locally  is a $p$-form,
with gauge invariance $A_p\to A_p+d\epsilon_{p-1}$ (with $\epsilon_{p-1}$
a $(p-1)$-form) and curvature $H=dA_p$.  But globally there can be 
non-trivial periods
$\int_D{H\over 2\pi}\in {\bf Z}$
for every $(p+1)$-cycle $D.$  More precisely,
$H$ is the de Rham representative of a characteristic class $x$ of $A_p$;
this class takes values in $H^{p+1}(M;{\bf Z})$ and can be an arbitrary
element of that group.
The Lagrangian, for a $p$-form field on
an $n$-manifold $M_n$, is
\begin{equation} 
I(H)={1\over 2\pi t}\int_{M_n}H\wedge *H,
\label{paction}
\end{equation}
with $t$ a positive constant.
In a more complete and rigorous description, the $A_p$ are ``differential
characters,'' for example $A_0$ is a map to ${\bf S}^1$, $A_1$ an
abelian gauge field, etc.  There is also a mathematical theory, not
yet much used by physicists, in which
a two-form field is understood as a connection on a gerbe, and the higher
$p$-forms are then related to more sophisticated objects.

We can compute the partition function as before.  We write $A_p=A_p'
+A_{p,h}$, where $A_p'$ is a globally defined $p$-form and $A_{p,h}$ is
a $p$-form field with harmonic curvature.  The curvature of $A_{p,h}$ is
determined by the characteristic class $x$ of $A_p$.
This leads to a description
of the partition function in which the interesting factor (for our purposes)
come from the sum over $x$.  It is\footnote{$\Theta$ is a function of the
metric on $M_n$, which enters through the induced metrix $\langle x,x
\rangle$ on the middle-dimensional cohomology.}
\begin{equation}
\Theta =
\sum_{x\in H^{p+1}(M_n;{\bf Z})}\exp\left(
-{\pi \over t}\langle x,x\rangle \right).
\label{htheta}
\end{equation}
As before $\langle x,x\rangle =\int_{M_n}x\wedge *x$.  The $*$ operator
that is used in this definition is only conformally invariant in the
middle dimension, so conformal invariance only holds if $n$ is even and
$p+1=n/2$.  Let us focus on this case.

If $n=4k$, then as we have already observed for $k=1$,  another term
${\theta\over (2\pi)^2}\int_{M_n}H\wedge H$
can be added to the action.  This leads to a modular function, similar to
what we have already described for $k=1$.

If $n=4k+2$, then ($H$ being a $(2k+1)$-form) $\int H\wedge H= 0$, so
we cannot add a $\theta$-term to the action.  But something else 
happens instead.

To understand this properly, we should at least temporarily return
to the case that $M_n$ is an $n$-manifold with {\it Lorentz} signature,
$-+++\dots +$, which is the real home of physics.  (In Lorentz signature
we normally
restrict $M_n$ to have a global Cauchy hypersurface, and no closed
timelike curves; normally, in Lorentz signature, we take $M_n$
 to have the
topology ${\bf R}\times M_{n-1}$, where ${\bf R}$ parametrizes the
``time'' and $M_{n-1}$ is ``space.'')
In $4k+2$ dimensions with Lorentz signature, a self-duality condition
$H=*H$ is possible for {\it real} $H$.  In $4k$ dimensions, self-duality
requires that $H$ be complex.  (In Euclidean signature, the conditions are reversed:
 a self-duality 
condition for a real middle-dimensional form is possible only in dimension
$4k$ rather than $4k+2$.  This result may be more familiar that the
corresponding Lorentzian statement.)

At any rate, in $4k+2$ dimensions with Lorentz signature, a middle-dimensional
classical $H$-field, obeying the Bianchi identity $dH=0$ and the 
Euler-Lagrange equation $d*H=0$, can be decomposed as $H=H_++H_-$, \pagebreak
where $H_\pm$ are real and
\begin{eqnarray}
 *H_\pm & = \pm H_\pm \nonumber \\
 dH_\pm & = 0.
\label{turnip}
\end{eqnarray}
Since classically it is consistent to set $H_-=0$, one may suspect that there
exists a quantum theory with $H_-=0$ and only $H_+$.  It turns out
that this is true if we choose the constant $t$ in the action eqn. \ref{paction}
properly.  

The lowest dimension of the form $4k+2$, to which this discussion is pertinent,
 is of course dimension two.
The self-dual quantum theory in dimension two has been extensively 
studied; it
is important in the
Segal-Frenkel-Kac vertex construction of representations of affine Lie
algebras, in bosonization of fermions and its applications to statistical
mechanics and representation theory, and in string theory.  In these
applications, it is important to consider generalizations of the  theory
we have considered to higher rank (by introducing several $H$ fields).
The generalization of picking a positive number $t$ is to pick a lattice with suitable
properties.  After dimension two, the next possibility (of the form $4k+2$)
is dimension six, and very interesting things, which we will indicate in section 3 below,
do occur in dimension six.   For understanding these phenomena, it is simplest and most useful
to set $t=1$.  However, theories with interesting (and in general more
complicated) properties can also be constructed for other rational
values of $t$.

There is no way to write a Lagrangian for the theory with $H_+$ only
-- since for example $\int_{M_{4k+2}}H_+\wedge H_+=0$.  This makes the quantum
theory subtle, but nevertheless it does exist, if we slightly relax our
axioms.  From the viewpoint that
we have been developing, this can be seen by writing the non-holomorphic
Siegel-Narain theta function of the lattice $\Lambda=H^{n/2}(M;{\bf Z})$, which
appears in eqn. \ref{htheta},
in terms of holomorphic theta functions.  For dimension $n=4k+2$, the lattice
$\Lambda$ has a {\it skew form} $(x,y)=\int x\wedge y$.  It also, of
course, just as in any other dimension, has the metric $\langle x,x\rangle
=\int_{M_{4n+2}}x\wedge *x$.  The skew form plus metric determine a complex
structure on the torus $T=H^{n/2}(M;U(1))/{\rm torsion}$.

Another important ingredient is a choice of ``quadratic refinement''
of the skew form.  A quadratic refinement of an integer-valued
 skew form $(x,y)$ is a ${\bf Z}_2$-valued
function $\phi:\Lambda\to{\bf Z}_2$ such that $\phi(x+y)=\phi(x)+\phi(y)
+(x,y)$ mod 2. There are $2^{b_{n/2}(M)}$ choices of such a $\phi$.
 Given a choice of $\phi$, by classical formulas one can 
construct a unitary line bundle with connection ${\cal L}_\phi\to T$ whose
curvature is the two-form determined by the skew form $(x,y)$.  This
turns $T$ into a ``principally polarized abelian variety,'' which has
an associated holomorphic theta function $\vartheta_\phi$.

It can be shown  (for a detailed discussion, see \cite{henning}) that the
non-holomorphic theta function $\Theta$ of eqn. \ref{htheta} which determines
the partition function of the original theory without self-duality
can be expressed in terms of the holomorphic theta functions 
$\vartheta_\phi$:
\begin{equation}
\Theta=\sum_\phi \,\vartheta_\phi \overline\vartheta_\phi.
\label{expresstheta}
\end{equation}
The sum runs over all choices of $\phi$.  If we could pick a $\phi$
in a natural way, we would interpret $\vartheta_\phi$ as the difficult
part, the
``numerator,'' of the partition function of the self-dual theory.
In fact, roughly speaking, a choice of a spin structure on $M$ determines
a $\phi$ (for more detail, see the last two papers in \cite{ewitten}, as well
as \cite{hopk} for an interpretation in terms of the Kervaire invariant).
So we modify the definition of conformal field theory to allow a choice
of spin structure and set the partition function $Z_{sd}$ of the self-dual theory to be
$Z_{sd}={\vartheta_\phi\over \det_+}.$
Here $\det_+$ is the result of projecting the determinant that comes
from the integral over topologically trivial fields onto the self-dual
part.  (Even in the absence of a self-dual projection, we did not discuss 
in any detail this determinant, which comes from the Gaussian integral
over the topologically trivial field $A_p'$.  For a discussion of it and
an explanation of its decomposition in self-dual and anti-self-dual factors
to get $\det_+$, see \cite{henning}.)

Many assertions we have made depend on
having set $t=1$.  For other values of $t$, to factorize $\Theta$ in terms of
holomorphic objects, we would need to use theta functions
at higher level; they would not be classified simply by a choice of
quadratic refinement; and the structure needed to pick a particular
holomorphic theta function would be more than a spin structure.

\bigskip\noindent{\it Relation Between $4k$ And $4k+2$ Dimensions}

My last goal in discussing these linear theories
 is to indicate, following \cite{verl}, how the existence of
a self-dual theory in $4k+2$ dimensions implies $SL(2,{\bf Z})$
(or $\Gamma_0(2)$) symmetry in $4k$ dimensions.  

First let us look at the situation classically. We formulate the 
$(4k+2)$-dimensional self-dual theory on the manifold $M_{4k+2}=M_{4k}\times
{\bf T}^2$, where $M_{4k}$ is a $(4k)$-manifold, and ${\bf T}^2$
a two-torus.  We take ${\bf T}^2={\bf R}^2/L$, where $L$ is a lattice
in the $u-v$ plane ${\bf R}^2$.  On ${\bf R}^2$ we take the metric
$ds^2= du^2+dv^2$.  So $E={\bf T}^2$ is an elliptic curve with a $\tau$
parameter $\tau_E$, which depends in the usual way on $L$.

Keeping the metric fixed on ${\bf T}^2$, we scale up the metric $g$ on 
$M_{4k}$
by $g\to \lambda g$, where we take $\lambda$ to become very large.
Any middle-dimensional form $H$ on $M_{4k}\times {\bf T}^2$ can be expanded
in Fourier modes on ${\bf T}^2$.  In our limit with ${\bf T}^2$ much smaller
than any characteristic radius of $M_{4k}$, the important modes
(which, for example, give the main contribution to the theta function) 
are constant, that is, invariant under translations on the torus.
So we can write
$H=F\wedge du+\tilde F\wedge dv + G + K\wedge du\wedge dv$
where $F,\tilde F, G, $ and $K$ are pullbacks from $M_{4k}$.

Self-duality of $H$ implies that $K=*G$ and that $\tilde F=*F$ (where
here $*$ is the duality operator on $M_{4k}$).
The $SL(2,{\bf Z})$ symmetry of ${\bf T}^2$ acts trivially on $G$ and $K$;
for that reason we have not much of interest to say about them.  Instead, we will
concentrate on $F$ and $\tilde F$.

The fact that $H$ is closed, $dH=0$, implies that $dF=d\tilde F=0$.
As $\tilde F=*F$, it follows that $dF=d*F=0$.  These are the usual conditions
(along with integrality of periods)
for $F$ to be the curvature of a $(2k-1)$-form field in $4k$ dimensions.
So, for example, if $k=1$, then $F$ is simply the curvature of an abelian
gauge field.  

So in the limit that the elliptic curve  $E$ is small compared to $M_{4k}$,
the self-dual theory on $M_{4k}\times E$, which I will call (a), is equivalent  to the
theory of a $(2k-1)$-form on $M_{4k}$ (plus less interesting contributions
from $G$ and $K$), which I will call (b).   

Suppose that this is true quantum mechanically.  The theory (a) depends
on the elliptic curve $E$, while (b) depends on $\tau =\theta/2\pi +4\pi
i/e^2$,
which modulo $SL(2,{\bf Z})$  determines an elliptic curve $E'$.

A natural guess is that $E\cong E'$, and if so (since theory (a) manifestly
depends only on $E$ and not on a contruction of $E$ using a specific
basis of the lattice $L$ or a specific $\tau$-parameter) this makes
obvious the $SL(2,{\bf Z})$ symmetry of theory (b).

The relation $E=E'$ can be established by comparing the theta functions.
But instead, I will motivate this relation in a way that will be helpful
when we study nonlinear theories in the next section.

Instead of reducing from $4k+2$ dimensions to $4k$ dimensions, let us
first compare $4k+2$ dimensions to $4k+1$ dimensions, and then take
a further step down to $4k$ dimensions.   So we formulate the self-dual
theory on $M_{4k+2}=M_{4k+1}\times {\bf S}^1$, with ${\bf S}^1$ described
by an angular variable $v$, $0\leq v\leq R$.  We fix the metric $dv^2$ on
${\bf S}^1$, and scale up the metric on $M_{4k+1}$ by a large factor.
In the limit, just as in the previous case,
we can assume $H=F\wedge dv+G$, where $F$ and $G$ are pullbacks from
$M_{4k+1}$.  Moreover, $G=*F$ and $dG=dF=0$, so $F$ obeys the conditions
$0=dF=*dF$ to be the curvature of an ``ordinary'' $(2k-1)$-form theory
on $M_{4k+1}$.\footnote{$2k-1$ is the degree of the potential, while the curvature $F$ is
of degree $2k$.}

Unlike the self-dual theory on $M_{4k+2}$, the ``ordinary'' theory on $M_{4k+1}$
does have a Lagrangian.  This Lagrangian depends on a free parameter
(called $t$ in eqn. \ref{paction}).  Conformal invariance on $M_{4k+1}\times {\bf S}^1$
implies that $t$ must be a constant multiple of $R$, so that the action
(apart from a constant that can be fixed by comparing the theta functions) is
\begin{equation}
I={1\over 4\pi R}\int_{M_{4k+1}}F\wedge *F.
\label{ilgo}
\end{equation}
The point of this formula is that if we rescale the metric of both
factors of $M_{4k+2}=M_{4k+1}\times {\bf S}^1$ by the same factor,
then $R$ (the circumference of ${\bf S}^1$) and $*$ (the Hodge $*$ operator
of $M_{4k+1}$ acting from $(2k)$-forms to $(2k+1)$-forms) scale in the same
way, so the action in eqn. \ref{ilgo} is invariant.   

The formula of eqn. \ref{ilgo} has the very unusual feature that $R$ is in the 
denominator.  If we had a Lagrangian in $4k+2$ dimensions, then after
specializing to $M_{4k+2}=M_{4k+1}\times {\bf S}^1$, we would deduce
what the action must be in $4k+1$ dimensions by simply ``integrating over
the fiber'' of the projection $M_{4k+2}\to M_{4k+1}$.  For fields that
are pullbacks from $M_{4k+1}$, this would inevitably give an action 
on $M_{4k+1}$ that is proportional to $R$ -- the volume of the fiber
-- and not to $R^{-1}$, as in eqn. \ref{ilgo}.  
But there is no classical action
in $4k+2$ dimensions, and the ``integration over the fiber'' is
a quantum operation that gives a factor of $R^{-1}$ instead of $R$.

Now let us return to the problem of comparing $4k+2$ dimensions to
$4k$ dimensions, and arguing that $E'$ is isomorphic to $E$.
We specialize to the case that the lattice $L$ is ``rectangular,''
generated  by the points $(S,0)$ and $(0,R)$
in the $u-v$ plane.  Accordingly, the torus $E\cong {\bf T}^2$ has
a decomposition as ${\bf S}\times {\bf S}'$, where ${\bf S}$ and ${\bf S}'$
are circles of circumference, respectively, $S$ and $R$.

We apply the previous reasoning to the decomposition
$M_{4k+2}=M_{4k+1}\times {\bf S}'$, with $M_{4k+1}=M_{4k}\times {\bf S}$.
Since ${\bf S}'$ has circumference $R$, the induced theory on $M_{4k+1}$
has action given by eqn. \ref{ilgo}.  Now, let us look at the decomposition
$M_{4k+1}=M_{4k}\times {\bf S}$.  Taking the length scale of $M_{4k}$ to
be large compared to that of ${\bf S}$, we would like to reduce to
a theory on $M_{4k}$.  For this step, since we do have a classical
action on $M_{4k+1}$, the reduction to a classical action on $M_{4k}$
is made simply by integrating over the fibers of the projection $M_{4k+1}\to
M_{4k}$.  As the fibers have volume $S$, the result is the
following action on $M_{4k}$:
\begin{equation}
I={1\over 4\pi}{S\over R}\int_{M_{4k}} F\wedge *F.
\label{naction}
\end{equation}
We see from eqn. \ref{naction} that the $\tau$ parameter of the theory on 
$M_{4k}$ is $\tau'=iS/R$.  But this in fact is the same as the $\tau$
parameter of the elliptic curve $E={\bf S}\times {\bf S}'$, so we have
demonstrated, for this example, that $E\cong E'$.  

In our two-step procedure of reducing from $M_{4k}\times {\bf S}\times 
{\bf S}'$, we made an arbitrary choice of reducing on ${\bf S}'$ first.
Had we proceeded in the opposite order, we would have arrived at 
$\tau'=iR/S$ instead of $iS/R$; the two results differ by the expected
modular transformation $\tau\to -1/\tau$.

One can extend the
above arguments to arbitrary $E$ with more work;
it is not necessary in this two-step reduction for
${\bf S}$ and ${\bf S}'$ to be orthogonal.  Of course, one can also make
the arguments more precise by study of the theta function of the
self-dual theory on $M_{4k+2}$. 

\section{Superconformal Field Theories In Four And Six Dimensions}

In $n$ dimensions, the conformal group of (conformally compactified)
Minkowski spacetime is
$SO(2,n)$.  A superconformal field theory, that is a conformal field
theory that is also supersymmetric, should have a supergroup $G$ of 
symmetries
whose bosonic part is $SO(2,n) \times K$, with $K$ a compact Lie group.
The fermionic part of the Lie algebra of $G$ 
should transform as a sum of spin representations of
$SO(2,n)$.  {\it A priori}, the spinors may appear in the Lie algebra
with any multiplicity, and for $n$ even, where $SO(2,n)$ has two distinct
spinor representations, these may appear with unequal multiplicities.

Nahm  considered the problem of classifying supergroups $G$ with
these properties.  The result is
 that solutions exist only for $n\leq 6$.  For $n=6$, the algebraic
 solution can be described as follows.  The group $G$ is $OSp(2,6|r)$ for
 some positive integer $r$.  Thus $K=Sp(r)$.  To describe the fermionic
 generators of $G$, first consider $G'=OSp(2,n|r)$ for general $n$.
 The fermionic generators of this group transform not as spinors
 but as the vector representation
 of $O(2,n)$ (tensored with the fundamental representation of $Sp(r)$).
 Thus for general $n$, the group $G'$ does not solve our algebraic 
 problem. However, precisely for $n=6$, we can use the {\it triality}
 symmetry of $O(2,6)$; by an outer automorphism of this group, its vector
 representation is equivalent to one of the two spinor representations.
 So modulo this automorphism, the group $G=OSp(2,6|r)$ does obey the
 right algebraic conditions and is a possible supergroup of symmetries
 for a superconformal field theory in six dimensions.

The algebraic solutions of Nahm's problem for $n<6$
 are similarly related to exceptional
isomorphisms  of Lie groups and supergroups of low rank.  (We give the example
of $n=4$ presently.)
Triality is in some sense the last of the exceptional isomorphisms, and the role of triality
for $n=6$ thus makes it plausible that $n=6$ is the maximum dimension
for superconformal symmetry, though I will not give a proof here.

As I  remarked in the introduction, this particular result by Nahm had little
immediate impact, since it was believed at the time that the correct
bound was really $n\leq 4$.
But in the mid-1990's, examples were found with $n=5,6$.
The known examples in dimension 6 have $r=1$ and $r=2$.
My goal in what follows will be to convey a few hints
about the $r=2$  examples.  A reference for some of what I will explain
is \cite{edward}.
 
\bigskip\noindent{\it Superconformal Gauge Theories In Four Dimensions}

We will need to know a few more facts about gauge theories in four dimensions.
The basic gauge theory with the standard Yang-Mills action
$I(A)={1\over 4e^2}\int \Tr F\wedge *F$ is conformally invariant
at the classical level, but not quantum mechanically.  There are many
ways to introduce additional fields and achieve quantum conformal 
invariance.

We will focus on superconformal field theories.  The superconformal
symmetries predicted by Nahm's analysis are $SU(2,2|{\cal N})$ for arbitrary
positive integer ${\cal N}$, as well
as an exceptional possibility $PSU(2,2|4)$. Note that $SU(2,2)$ is isomorphic
to $SO(2,4)$, and that the fermionic part of the super Lie algebra of
$SU(2,2|{\cal N})$ (or of $PSU(2,2|4)$) transforms as ${\cal N}$ 
copies of $V\oplus
\overline V$, where $V$ is the defining four-dimensional representation
of $SU(2,2)$.  $V$ and $\overline V$ are isomorphic to the two spinor 
representations of $SO(2,4)$, so $SU(2,2|{\cal N})$ and $PSU(2,2|4)$ do
solve the algebraic problem posed by Nahm.  The supergroups 
$SU(p,q|{\cal N})$
exist for all positive integers $p,q,{\cal N}$, but it takes the exceptional
isomorphism $SU(2,2)\cong SO(2,4)$ to get a solution of the problem
considered by Nahm.

Examples of superconformal field theories in four dimensions exist for ${\cal 
N}=1,2,$ and 4.  For ${\cal N}=1$, there are myriads of
possibilities -- though much more constrained than in the absence 
of supersymmetry -- while the examples with  ${\cal N}=2$ and 
${\cal N}=4$ are so highly constrained that a complete classification is 
possible. In particular, for ${\cal N}=4$, the fields that must be 
included are completely determined by the choice of the gauge 
group $G$.  For ${\cal N}=2$, one also picks a representation of $G$ that 
obeys a certain condition on the trace of the quadratic Casimir 
operator (there are finitely many choices for each given $G$).  
We will concentrate on the examples with ${\cal N}=4$; they 
have the exceptional  $PSU(2,2|4)$ symmetry.

\bigskip\noindent{\it ${\cal N}=4$ Super Yang-Mills Theory}

The fields of ${\cal N}=4$ super Yang-Mills theory are the gauge field
$A$ plus fermion and scalar fields required by the supersymmetry.
The Lagrangian is
\begin{equation}
I(A,\dots)=\int_{M_4}\Tr\left({1\over 4e^2}F\wedge *F+{i\theta\over
8\pi^2}F\wedge F+\dots\right).
\label{gffs}
\end{equation}
where the ellipses refer to terms involving the additional fields. 
 
If we set $\tau={\theta\over 2\pi}+{4\pi i\over e^2}$,
then the Montonen-Olive duality conjecture \cite{monto} asserts an $SL(2,{\bf Z})$
symmetry acting on $\tau$.  Actually, the element 
\begin{equation}
\left( \begin{array}{cc}
0 & 1\\
-1 & 0 \label{elem} \end{array} \right)
\end{equation}
of $SL(2,{\bf Z})$ is conjectured to map the ${\cal N}=4$ theory with
gauge group $G$ to the same theory with the Langlands dual group, while
also mapping $\tau$ to $-1/\tau$.
So in general the precise modular properties are a little involved,
somewhat analogous to the fact that in section 2, we found in general
$\Gamma_0(2)$ rather than full $SL(2,{\bf Z})$ symmetry.
By around 1995, many developments in the study of supersymmetric gauge
theories and string theories gave strong support for the Montonen-Olive
conjecture.

If we formulate the ${\cal N}=4$ theory on a compact four-manifold $M$,
endowed with some metric tensor $g$, the partition function $Z(M,g;\tau)$
is, according to the Montonen-Olive conjecture, a modular function of
$\tau$.  It is not in general holomorphic or anti-holomorphic in $\tau$,
and it depends non-trivially on $g$, so it is not a topological invariant
of $M$.

However \cite{vafawitt}, there is a ``twisted'' version of the theory that is a 
topological
field theory and still $SL(2,{\bf Z})$-invariant.  For a four-manifold
$M$ with $b_{2,+}(M)>1$, the partition function is holomorphic (with a pole
at the ``cusp'') and a topological invariant of $M$.
In fact, setting $q=\exp(2\pi i\tau)$, the partition function can be written
\begin{equation}
Z(M;\tau)=q^{-c}\sum_{n=0}^\infty a_nq^n,
\label{umbo}
\end{equation}
where, assuming a certain vanishing theorem holds, $a_n$ is the Euler
characteristic of the moduli space of $G$-instantons of instanton
number $n$.  In general, $a_n$ is the ``number'' of solutions, weighted
by sign, for a certain coupled system of equations for the connection
plus certain additional fields.  These more elaborate equations, which
are somewhat analogous to the Seiberg-Witten equations and
have similarly nice Bochner formulas (related in both cases to supersymmetry),
were described in \cite{vafawitt}.

\bigskip\noindent{\it Explanation From Six Dimensions}

So if the $SL(2,{\bf Z})$ conjecture of Montonen and Olive holds, 
the functions defined in eqn. \ref{umbo} are modular.  But why should the
${\cal N}=4$ supersymetric gauge theory in four dimensions have $SL(2,{\bf
Z})$ symmetry?

Several explanations emerged from string theory work in the mid-1990's.
Of these, one \cite{edward} is in the spirit of what we discussed for linear theories
in section 2.  In its original form, this explanation only works for
simply-laced $G$, that is for $G$ of type $A,D$, or $E$.  I will limit
the following discussion to this case.  (For simply-laced $G$, $G$ is
locally isomorphic to its Langlands dual, and the statement of Montonen-Olive
duality becomes simpler.)

The surprise which leads to an insight about Montonen-Olive duality is that
in dimension $n=6$, there is for each choice of simply-laced group $G$
a superconformal field theory that is a sort of nonlinear (and supersymmetric)
version of the self-dual theory that we discussed in section 2.  This
exotic six-dimensional theory was found originally \cite{edward} by considering
Type IIB superstring theory at an $A-D-E$ singularity.  

The superconformal symmetry of this theory is the supergroup 
$OSp(2,6|2)$.  When it
is formulated on a six-manifold $M_6=M_4\times E$, with $E$ an elliptic
curve, the resulting behavior is quite similar to what we have discussed
in section 2 for the linear self-dual theory.  Taking a product metric
on $M_4\times E$, in the limit that $M_4$ is much larger than $E$,
the six-dimensional theory reduces to the four-dimensional ${\cal N}=4$
theory with gauge group $G$ and $\tau$ parameter determined by $E$.  
Just as in section 2, this makes manifest the Montonen-Olive symmetry
of the ${\cal N}=4$ theory.  From this point of view, Montonen-Olive symmetry
 reflects the fact that the six-dimensional theory on $M_4\times E$ depends
 only on $E$ and not on a specific way of constructing $E$ using 
 a $\tau $ parameter.
 
Further extending the analogy with what we discussed in section two for
linear theories, if we formulate this theory on $M_5\times {\bf S}$, where
${\bf S}$ is a circle of circumference $R$, we get at  distances large
compared to $R$  a five-dimensional gauge theory, with gauge group
$G$, and action proportional to $R^{-1}$ rather than $R$.  
As in section 2, this shows
that the five-dimensional action cannot be obtained by a classical
process of ``integrating over the fiber''; it gives an
obstruction to deriving the six-dimensional theory from a Lagrangian.

 The six-dimensional theory that comes from Type IIB superstring theory
 at the $A-D-E$ singularity might be called a ``nonabelian gerbe theory,''
 as it is an analog for  $A-D-E$ groups of the linear theory
 discussed in section two with a two-form field and a self-dual three-form
 curvature.  Under a certain perturbation (to a vacuum with spontaneous
 symmetry breaking in six dimensions), the six-dimensional $A-D-E$
 theory  reduces at low energies to a theory that can be described more
 explicitly; this theory is a  more elaborate version of the theory with
 self-dual curvature that
 we considered in section 3.  In this theory, the gerbe-like field
 has a characteristic class that
  takes values not in $H^3(M;{\bf Z})$, but in $H^3(M;
 {\bf Z})\otimes \Lambda$, where $\Lambda$ is the root lattice of $G$.
 Physicists describe this roughly by saying that, if $r$ denotes the rank
 of $G$, there are $r$ self-dual two-form fields (i.e., two-form fields
 whose curvature is a self-dual three-form).  
 
The basic hallmark of the six-dimensional theory is that on the one
hand it can be perturbed to give something that  we recognize as a
gerbe theory of rank $r$; on the other hand, it can be perturbed to give
non-abelian gauge theory with gauge group $G$.  Combining the two facts,
this six-dimensional theory is a sort of quantum nonabelian gerbe theory.
I doubt very much that this structure is accessible in the world of
classical geometry; it belongs to the realm of quantum field theory.
But it has manifestations in the classical world, such as the modular
nature of the generating function (eqn. \ref{umbo}) 
of Euler characteristics of
instanton moduli spaces.


\begin{thebibliography}{99}


\bibitem{gsegal}
G. Segal, ``The Definition Of Conformal Field Theory,'' in
{\it Differential Geometric Methods In Theoretical Physics}
(proceedings, Como, 1987), and these proceedings.

\bibitem{wnahm}
W. Nahm, ``Supersymmetries And Their Representations,'' Nucl. Phys.
{\bf B135} (1978) 149.

\bibitem{ewitten}
E. Witten, ``On $S$-Duality In Abelian Gauge Theory,'' hep-th/9505186, Selecta Mathematica
{\bf 1} (1995) 383;
``Five-Brane Effective Action In $M$-Theory,'' hep-th/9610234, J. Geom.
Phys. {\bf 22} (1997) 103; ``Duality Relations Among
Topological Effects In String Theory,''
hep-th/9912086, JHEP {\bf 0005:031,2000}.

\bibitem{henning}
M. Henningson, ``The Quantum Hilbert Space Of
A Chiral Two Form In $D=(5+1)$ Dimensions, hep-th/0111150,
JHEP {\bf 0203:021,2002}.

\bibitem{hopk}
M. Hopkins and I. M. Singer, ``Quadratic Functions In Geometry,
Topology, and $M$-Theory,'' math.at/0211216.

\bibitem{verl}
E. Verlinde, ``Global Aspects Of Electric-Magnetic
Duality,'' hep-th/9506011, Nucl. Phys. {\bf B455} (1995) 211.

\bibitem{edward}
E. Witten, ``Some Comments On String Dynamics,'' hep-th/9507121, in
{\it Future Perspectives In String Theory}, (World Scientific, 1996),
ed. I. Bars et. al., 501.

\bibitem{monto}
C. Montonen and D. Olive, ``Magnetic Monopoles As Gauge Particles?''
Phys. Lett. {\bf B72} (1977) 117; P. Goddard, J. Nuyts, and D. Olive,
``Gauge Theories And Magnetic Charge,'' Nucl. Phys. {\bf B125} (1977) 1;
H. Osborn, ``Topological Charges For ${\cal N}=4$ Supersymmetric
Gauge Theories And Monopoles Of Spin 1,'' Phys. Lett. {\bf B83} (1979) 321.

\bibitem{vafawitt}
C. Vafa and E. Witten, ``A Strong Coupling Test Of $S$-Duality,''
hep-th/9408074, Nucl. Phys. {\bf B431} (1994) 3.

\end{thebibliography}
\end{document}